\documentstyle[amscd,amssymb]{amsart}
\newcommand{\nc}{\newcommand}

\nc{\ad}{{\mbox{\bf{ad}}}}
\nc{\AJ}{{\operatorname{aj}}}
\nc{\Aut}{{\operatorname{Aut}}}
\nc{\Bls}{{{\cal B}ls}}
\nc{\Boxtimes}{{\fbox{$\times$}}}
\nc{\blt}{{\bullet}}
\nc{\bSt}{{\mbox{\bf{St}}}}
\nc{\card}{{\operatorname{card}}}
\nc{\Cch}{{\check{C}}}
\nc{\cd}{{\operatorname{cd}}}
\nc{\Ch}{{\operatorname{Ch}}}
\nc{\chara}{{\operatorname{char}}}
\nc{\CHom}{{\cal{H}om}}
\nc{\Coker}{{\operatorname{Coker}}}
\nc{\codim}{{\operatorname{codim}}}
\nc{\Cone}{{\operatorname{Cone}}}
\nc{\cSgn}{{\cal{S}gn}}
\nc{\depth}{{\operatorname{depth}}}
\nc{\dirlim}{{\underset{\rightarrow}{\operatorname{lim}}}}
\nc{\dotbox}{{\overset{\bullet}{\boxtimes}}}
\nc{\dotimes}{{\overset{\bullet}{\otimes}}}
\nc{\Ed}{{\operatorname{Edge}}}
\nc{\emp}{{\emptyset}}
\nc{\Ext}{{\operatorname{Ext}}}
\nc{\Fac}{{\cal{F}ac}}
\nc{\Fun}{{\operatorname{F}}}
\nc{\FS}{{\cal{FS}}}
\nc{\Hom}{{\operatorname{Hom}}}
\nc{\had}{{{\hat{\mbox{\bf{ad}}}}}}
\nc{\hgt}{{\operatorname{ht}}}
\nc{\Id}{{\operatorname{Id}}}
\nc{\id}{{\operatorname{id}}}
\nc{\Ima}{{\operatorname{Im}}}
\nc{\ind}{{\operatorname{ind}}}
\nc{\Ind}{{\operatorname{Ind}}}
\nc{\infi}{{\operatorname{inf}}}
\nc{\infh}{{\frac{\infty}{2}}}
\nc{\invlim}{{\underset{\leftarrow}{\operatorname{lim}}}}
\nc{\Jac}{{{\cal J}ac}}
\nc{\Ker}{{\operatorname{Ker}}}
\nc{\lcm}{{\operatorname{lcm}}}
\nc{\length}{{\operatorname{length}}}
\nc{\Locsys}{{{\cal L}ocsys}}
\nc{\Map}{{{\cal M}ap}}
\nc{\modul}{{\operatorname{mod}}}
\nc{\Mor}{{\operatorname{Mor}}}
\nc{\MS}{{\cal{MS}}}
\nc{\Ob}{{\operatorname{Ob}}}
\nc{\opp}{{\operatorname{opp}}}
\nc{\Or}{{{\cal O}r}}
\nc{\Ord}{{{\cal O}rd}}
\nc{\Part}{{{\cal P}art}}
\nc{\PGL}{{\operatorname{PGL}}}
\nc{\Pic}{{\operatorname{Pic}}}
\nc{\Rep}{{{\cal{R}}ep}}
\nc{\rk}{{\operatorname{rk}}}
\nc{\Sets}{{{\cal{S}}ets}}
\nc{\Sew}{{{\cal{S}}ew}}
\nc{\sgn}{{\operatorname{sgn}}}
\nc{\Sh}{{{\cal S}h}}
\nc{\Sign}{{{\cal S}ign}}
\nc{\Spe}{{\mbox{\bf{Sp}}}}
\nc{\supr}{{\operatorname{sup}}}
\nc{\Supp}{{\operatorname{Supp}}}
\nc{\supp}{{\operatorname{supp}}}
\nc{\Teich}{{{\cal{T}}eich}}
\nc{\tFS}{{\widetilde{\cal{FS}}}}
\nc{\Tor}{{\operatorname{Tor}}}
\nc{\totimes}{{\tilde{\otimes}}}
\nc{\tr}{{\operatorname{tr}}}
\nc{\tRep}{{\widetilde{{\cal R}ep}}}
\nc{\tTeich}{{\widetilde{{\cal T}eich}}}
\nc{\Vc}{{\mbox{\bf{V}}_{\mbox{\bf{c}}}}}
\nc{\Vect}{{{\cal V}ect}}
\nc{\Ve}{{\operatorname{Vert}}}
\nc{\wt}{{\widetilde}}

\nc{\bo}{{\mbox{\bf{0}}}}
\nc{\One}{{\mbox{\bf{1}}}}
\nc{\one}{{\mbox{\bf{1}}}}

\nc{\BA}{{\Bbb A}}
\nc{\bA}{{\overline{A}}}
\nc{\ba}{{\mbox{\bf{a}}}}
\nc{\baB}{{\overline{B}}}
\nc{\baeta}{{\bar{\eta}}}
\nc{\baJ}{{\bar{J}}}
\nc{\BB}{{\Bbb B}}
\nc{\bB}{{\mbox{\bf{B}}}}
\nc{\bc}{{\mbox{\bf{c}}}}
\nc{\bC}{{\overline{C}}}
\nc{\BC}{{\Bbb{C}}}
\nc{\bCC}{{\overline{\cal{C}}}}
\nc{\bCM}{{\overline{\cal{M}}}}
\nc{\bD}{{\bar{D}}}
\nc{\BD}{{\overline{D}}}
\nc{\bd}{{\mbox{\bf{d}}}}
\nc{\BE}{{\overline{E}}}
\nc{\BF}{{\overline{F}}}
\nc{\bF}{{\mbox{\bf{F}}}}
\nc{\bg}{{\mbox{\bf{g}}}}
\nc{\bG}{{\mbox{\bf{G}}}}
\nc{\BG}{{\Bbb G}}
\nc{\bGamma}{{\overline{\Gamma}}}
\nc{\bbH}{{     {\mbox{\bf{H}}}_a       }}
\nc{\bH}{{\mbox{\bf{H}}}}
\nc{\bI}{{\mbox{\bf{I}}}}
\nc{\bL}{{\mbox{\bf{L}}}}
\nc{\BL}{{\Bbb{L}}}
\nc{\blambda}{{\bar{\lambda}}}
\nc{\bM}{{\mbox{\bf{V}}}}
\nc{\bmu}{{\vec{\mu}}}
\nc{\bN}{{\mbox{\bf{N}}}}
\nc{\BN}{{\Bbb{N}}}
\nc{\bnu}{{\mbox{\boldmath{${\nu}$}}}}
\nc{\bof}{{\mbox{\bf{f}}}}
\nc{\BP}{{\Bbb P}}
\nc{\bP}{{\mbox{\bf{P}}}}
\nc{\BPO}{{\overset{\circ}{\BP}}}
\nc{\BQ}{{\Bbb Q}}
\nc{\BR}{{\Bbb{R}}}
\nc{\bR}{{\mbox{\bf{R}}}}
\nc{\bp}{{\mbox{\bf{p}}}}
\nc{\barq}{{\mbox{\bf{q}}}}
\nc{\br}{{\mbox{\bf{r}}}}
\nc{\breta}{{\bar{\eta}}}
\nc{\bs}{{\mbox{\bf{s}}}}
\nc{\bS}{{\mbox{\bf{S}}}}
\nc{\bt}{{\mbox{\bf{t}}}}
\nc{\BT}{{\Bbb T}}
\nc{\bU}{{\mbox{\bf{U}}}}
\nc{\bV}{{\mbox{\bf{V}}}}
\nc{\bu}{{\mbox{\bf{u}}}}
\nc{\BUpsilon}{{\bar{\Upsilon}}}
\nc{\bw}{{\mbox{\bf{w}}}}
\nc{\bx}{{\mbox{\bf{x}}}}
\nc{\bX}{{\mbox{\bf{X}}}}
\nc{\BZ}{{\Bbb{Z}}}
\nc{\bz}{{\mbox{\bf{z}}}}
\nc{\bZ}{{\mbox{\bf{Z}}}}
\nc{\bzero}{\mbox{\boldmath{$0$}}}

\nc{\CA}{{\cal A}}
\nc{\CAD}{{\overset{\bullet}{\cal{A}}}}
\nc{\CAO}{{\overset{\circ}{\cal{A}}}}
\nc{\CB}{{\cal B}}
\nc{\CC}{{\cal C}}
\nc{\CalD}{{\cal D}}
\nc{\CE}{{\cal E}}
\nc{\CF}{{\cal F}}
\nc{\CG}{{\cal G}}
\nc{\CH}{{\cal H}}
\nc{\CI}{{\cal I}}
\nc{\CID}{{\overset{\bullet}{\cal{I}}}}
\nc{\CJ}{{\cal J}}
\nc{\CK}{{\cal K}}
\nc{\CL}{{\cal L}}
\nc{\CM}{{\cal M}}
\nc{\CN}{{\cal N}}
\nc{\CO}{{\cal O}}
\nc{\CP}{{\Bbb P}_\gamma}
\nc{\CPO}{{\overset{\circ}{\cal{P}}}}
\nc{\CQ}{{\cal Q}}
\nc{\CR}{{\cal R}}
\nc{\CS}{{\cal S}}
\nc{\CT}{{\cal T}}
\nc{\CTD}{{\overset{\bullet}{\cal{T}}}}
\nc{\CTPO}{{\overset{\circ}{\cal{T}\cal{P}}}}
\nc{\CU}{{\cal{U}}}
\nc{\CV}{{\cal V}}
\nc{\CW}{{\cal W}}
\nc{\CX}{{\cal X}}
\nc{\CY}{{\cal Y}}
\nc{\CZ}{{\cal Z}}

\nc{\dCL}{{\overset{\bullet}{\cal{L}}}}
\nc{\dd}{{\operatorname{d}}}
\nc{\ddelta}{{\overset{\bullet}{\delta}}}
\nc{\dfu}{{\overset{\bullet}{\frak{u}}}}
\nc{\dlambda}{{\overset{\bullet}{\lambda}}}
\nc{\DO}{{\overset{\circ}{D}}}
\nc{\dpar}{{\partial}}
\nc{\dS}{{\overset{\bullet}{S}}}
\nc{\dT}{{\overset{\bullet}{T}}}

\nc{\fA}{{\frak{A}}}
\nc{\fb}{{\frak{b}}}
\nc{\fC}{{\frak{C}}}
\nc{\fD}{{\frak{D}}}
\nc{\fd}{{\frak{d}}}
\nc{\fE}{{\frak{E}}}
\nc{\fF}{{\frak{F}}}
\nc{\fEF}{{\frak{EF}}}
\nc{\fFE}{{\frak{FE}}}
\nc{\ff}{{\frak{f}}}
\nc{\fg}{{\frak{g}}}
\nc{\fG}{{\frak{G}}}
\nc{\fH}{{\frak{H}}}
\nc{\fl}{{\frak{l}}}
\nc{\fK}{{\frak{K}}}
\nc{\fL}{{\frak{L}}}
\nc{\fM}{{\frak{M}}}
\nc{\fN}{{\frak{N}}}
\nc{\fn}{{\frak{n}}}
\nc{\fp}{{\frak{p}}}
\nc{\fu}{{\frak{u}}}
\nc{\fZ}{{\frak{Z}}}

\nc{\hCH}{{\hat{\cal{H}}}}
\nc{\hCI}{{\hat{\cal{I}}}}
\nc{\hfC}{{\hat{\frak{C}}}}
\nc{\hfg}{{\hat{\frak{g}}}}
\nc{\hL}{{\hat{L}}}
\nc{\HO}{{\overset{\circ}{H}}}
\nc{\hpsi}{{\hat{\psi}}}
\nc{\hx}{{\hat{x}}}

\nc{\jo}{{\overset{\circ}{j}}}

\nc{\phid}{{\overset{\bullet}{\phi}}}

\nc{\tA}{{\tilde{A}}}
\nc{\ta}{{\tilde{a}}}
\nc{\tB}{{\tilde{B}}}
\nc{\tb}{{\tilde{b}}}
\nc{\tBP}{{\tilde{\BP}}}
\nc{\tC}{{\tilde{C}}}
\nc{\tc}{{\tilde{c}}}
\nc{\tCA}{{\tilde{\cal{A}}}}
\nc{\tCC}{{\tilde{\cal{C}}}}
\nc{\tCH}{{\tilde{\cal{H}}}}
\nc{\tCI}{{\tilde{\cal{I}}}}
\nc{\tCO}{{\tilde{\cal{O}}}}
\nc{\tCP}{{\tilde{\cal{P}}}}
\nc{\tCT}{{\tilde{\cal{T}}}}
\nc{\tD}{{\tilde{D}}}
\nc{\tDelta}{{\tilde{\Delta}}}
\nc{\tE}{{\tilde E}}
\nc{\tF}{{\tilde F}}
\nc{\tfD}{{\tilde{\frak{D}}}}
\nc{\tfF}{{\tilde{\frak{F}}}}
\nc{\tff}{{\tilde{\frak{f}}}}
\nc{\tfu}{{\tilde{\frak{u}}}}
\nc{\tJ}{{\tilde{J}}}
\nc{\tj}{{\tilde{j}}}
\nc{\tK}{{\tilde K}}
\nc{\tL}{{\tilde{L}}}
\nc{\tM}{{\tilde{M}}}
\nc{\tP}{{\tilde{P}}}
\nc{\tPhi}{{\tilde{\Phi}}}
\nc{\tpi}{\tilde{\pi}}
\nc{\TPO}{{\overset{\circ}{T\BP}}}
\nc{\tR}{{\tilde{R}}}
\nc{\tS}{{\tilde S}}
\nc{\tT}{{\tilde{T}}}
\nc{\ttau}{{\tilde{\tau}}}
\nc{\ttheta}{{\tilde{\theta}}}
\nc{\tU}{{\tilde{U}}}
\nc{\tUpsilon}{{\tilde{\Upsilon}}}
\nc{\tW}{{\tilde W}}
\nc{\ty}{{\tilde y}}
\nc{\tY}{{\tilde Y}}
\nc{\txi}{{\tilde{\xi}}}

\nc{\UD}{{\overset{\bullet}{U}}}
\nc{\UO}{{\overset{\circ}{U}}}
\nc{\uQ}{{\underline{\Bbb Q}}}

\nc{\vA}{{\vec{A}}}
\nc{\valpha}{{\vec{\alpha}}}
\nc{\vbeta}{{\vec{\beta}}}
\nc{\vc}{{\vec{c}}}
\nc{\vD}{{\vec{D}}}
\nc{\vd}{{\vec{d}}}
\nc{\vgamma}{{\vec{\gamma}}}
\nc{\vK}{{\vec{K}}}
\nc{\vlambda}{{\vec{\lambda}}}
\nc{\vmu}{{\vec{\mu}}}
\nc{\vnu}{{\vec{\nu}}}
\nc{\vo}{{\vec{0}}}
\nc{\vu}{{\vec{u}}}
\nc{\vx}{{\vec{x}}}
\nc{\vy}{\vec{y}}
\nc{\vzero}{\vec{0}}

\nc{\XO}{{\overset{\circ}{X}}}

\nc{\ya}{{\operatorname{aj}}}

\nc{\nen}{\newenvironment}
\nc{\ol}{\overline}
\nc{\ul}{\underline}
\nc{\ra}{\rightarrow}
\nc{\lra}{\longrightarrow}
\nc{\Lra}{\Longrightarrow}
\nc{\lla}{\longleftarrow}
\nc{\Llra}{\Longleftrightarrow}
\nc{\hra}{\hookrightarrow}
\nc{\iso}{\overset{\sim}{\lra}}
\nc{\rlh}{\rightleftharpoons}

\nc{\IC}{{\cal{IC}}}
\nc{\PS}{{\cal{PS}}}
\nc{\oCQ}{{\overline{\cal Q}}}
\nc{\oCZ}{{\overline{\cal Z}}}
\nc{\dZ}{{\overset{\bullet}{\cal Z}}{}}
\nc{\oZ}{{\overset{\circ}{\cal Z}}{}}
\nc{\dP}{{\overset{\bullet}{\cal P}}{}}
\nc{\oP}{{\overset{\circ}{\cal P}}{}}
\nc{\Ue}{{U_\varepsilon}}
\nc{\Upe}{{\Upsilon_\varepsilon}}
\nc{\crho}{{\check{\rho}}}
\nc{\ctheta}{{\check{\theta}}}

\nc{\QD}[1]{\CQ^D_{#1}}
\nc{\QL}[1]{\CQ^L_{#1}}
\nc{\QK}[1]{\CQ^K_{#1}}
\nc{\al}{\alpha}
\nc{\ga}{\gamma}
\nc{\ka}{\kappa}
\nc{\pr}{\bp\times \br}
\nc{\ME}{{\check \fE}}
\nc{\PP}{{\Bbb P}}
\nc{\CCC}{{\Bbb C}}
\nc{\NNN}{{\Bbb N}}
\nc{\ZZZ}{{\Bbb Z}}
\nc{\RT}{{\frak R\frak T}}
\nc{\Gr}{\operatorname{\text{Gr}}}
\nc{\wti}{\widetilde}
\nc{\wha}{\widehat}
\nc{\vphi}{\varphi}
\nc{\deff}{\operatorname{\text{def}}}
\nc{\hnu}{{\hat\nu}}
\nc{\tnu}{{\tilde\nu}}
\nc{\tka}{{\tilde\ka}}
\nc{\ti}{\tilde}
\nc{\fS}{{\frak S}}
\nc{\EGo}{{\overset{\circ}{\fE^\al_\Gamma}}}
\nc{\eac}{\fE^\al_{\underbrace{\theta_1,\dots,\theta_1}_{c_1},\dots,
\underbrace{\theta_\nu,\dots,\theta_\nu}_{c_\nu}}}
\nc{\lbr}{\{\!\{}
\nc{\rbr}{\}\!\}}

\nc{\Thm}[1]{Theorem~\ref{#1}}
\nc{\Prop}[1]{Proposition~\ref{#1}}
\nc{\Lem}[1]{Lemma~\ref{#1}}
\nc{\Cor}[1]{Corollary~\ref{#1}}
\nc{\Conj}[1]{Conjecture~\ref{#1}}
\nc{\Claim}[1]{Claim~\ref{#1}}
\nc{\Defn}[1]{Definition~\ref{#1}}
\nc{\Exa}[1]{Example~\ref{#1}}
\nc{\Rem}[1]{Remark~\ref{#1}}
\nc{\Note}[1]{Note~\ref{#1}}

\nen{thm}[1]{\label{#1}{\bf Theorem.\ } \em}{}
\nen{prop}[1]{\label{#1}{\bf Proposition.\ } \em}{}
\nen{lem}[1]{\label{#1}{\bf Lemma.\ } \em}{}
\nen{cor}[1]{\label{#1}{\bf Corollary.\ } \em}{}
\nen{conj}[1]{\label{#1}{\bf Conjecture.\ } \em}{}
\nen{claim}[1]{\label{#1}{\bf Claim.\ } \em}{}

\nen{defn}[1]{\label{#1}{\bf Definition.\ } }{}
\nen{exa}[1]{\label{#1}{\bf Example.\ } }{}

\nen{rem}[1]{\label{#1}{\em Remark.\ } }{}
\nen{note}[1]{\label{#1}{\em Note.\ } }{}
\nen{exer}[1]{\label{#1}{\em Exercise.\ } }{}

\nc{\TE}{{\tilde E}}
\nc{\opl}{\oplus}
\nc{\opll}{\mathop{\opl}\limits}
\nc{\ot}{\otimes}
\nc{\Om}{\Omega}

\begin{document}

\title[]{The Singular Supports of IC sheaves on Quasimaps' Spaces are
irreducible}
\author{Michael Finkelberg}
\address{Independent Moscow University, 11 Bolshoj Vlasjevskij pereulok,
Moscow 121002 Russia}
\email{fnklberg@@main.mccme.rssi.ru}
\author{Alexander Kuznetsov}
\address{Independent Moscow University, 11 Bolshoj Vlasjevskij pereulok,
Moscow 121002 Russia}
\email{sasha@@ium.ips.ras.ru}
\author{Ivan Mirkovi\'c}
\address{Dept. of Mathematics and Statistics, University of Massachusetts at
Amherst, Amherst MA 01003-4515, USA}
\email{mirkovic@@math.umass.edu}
\thanks{M.F. and A.K. were partially supported by CRDF grant RM1-265.
M.F.~was partially supported by INTAS-94-4720. I.M. was partially supported
by NSF}
\maketitle

\section{Introduction}

\subsection{}
Let $C$ be a smooth projective curve of genus 0. Let $\CB$ be the variety
of complete flags in an $n$-dimensional vector space $V$.
Given an $(n-1)$-tuple $\alpha\in\BN[I]$
of positive integers one can consider the space $\CQ_\alpha$ of algebraic
maps of degree $\alpha$ from $C$ to $\CB$. This space is noncompact. Some
remarkable compactifications $\CQ^D_\alpha$ (Quasimaps),
$\CQ^L_\alpha$ (Quasiflags) of
$\CQ_\alpha$ were constructed by Drinfeld and Laumon respectively.
In ~\cite{k} it was proved that the natural map $\pi:\ \CQ^L_\alpha\to
\CQ^D_\alpha$ is a small resolution of singularities. The aim of the present
note is to study the singular support of the Goresky-MacPherson sheaf
$IC_\alpha$ on the Quasimaps' space $\CQ^D_\alpha$.

Namely, we prove that this singular support $SS(IC_\alpha)$ is irreducible. The
proof
is based on the {\em factorization property} of Quasimaps' space and on the
detailed analysis of Laumon's resolution $\pi:\ \CQ^L_\alpha\to\CQ^D_\alpha$.

We are grateful to P.Schapira for the illuminating correspondence.

This note is a sequel to ~\cite{k} and ~\cite{fk}.
In fact, the local geometry of $\CQ^D_\alpha$ was the subject of ~\cite{k};
the global geometry of $\CQ^D_\alpha$ was the subject of ~\cite{fk},
while the microlocal geometry of $\CQ^D_\alpha$ is the subject of the present
work. We will freely refer the reader to ~\cite{k} and ~\cite{fk}.

\section{Reductions of the main theorem}

\subsection{Notations}

\subsubsection{}
\label{not}
We choose a basis $\{v_1,\ldots,v_n\}$ in $V$. This choice
defines a Cartan subgroup $H\subset G=SL(V)=SL_n$ of matrices diagonal with
respect to
this basis, and a Borel subgroup $B\subset G$ of matrices upper triangular
with respect to this basis. We have $\CB=G/B$.

Let $I=\{1,\ldots,n-1\}$ be the set of simple coroots of $G=SL_n$.
Let $R^+$ denote the set of positive coroots,
and let $2\rho=\sum_{\theta\in R^+}\theta$.
For $\alpha=\sum a_ii\in\BN[I]$ we set $|\alpha|:=\sum a_i$.
Let $X$ be the lattice of weights of $G,H$. Let $X^+\subset X$ be the set of
dominant (with respect to $B$) weights. For $\lambda\in X^+$ let $V_\lambda$
denote the irreducible representation of $G$ with the highest weight $\lambda$.

Recall the notations of ~\cite{k} concerning Kostant's partition function.
For $\gamma\in\BN[I]$ a {\em Kostant partition} of $\gamma$ is a decomposition
of $\gamma$ into a sum of positive coroots with multiplicities.
The set of Kostant partitions of $\gamma$ is denoted by
$\fK(\gamma)$.

There is a natural bijection between the set of pairs $1\leq q\leq p\leq n-1$
and $R^+$, namely, $(p,q)$ corresponds to $i_q+i_{q+1}+\ldots+i_p$. Thus a
Kostant partition $\kappa$ is given by a collection of nonnegative integers
$(\kappa_{p,q}), 1\leq q\leq p\leq n-1$.
Following {\em loc. cit.} (9) we define a collection $\mu(\kappa)$ as follows:
$\mu_{p,q}=\sum_{r\leq q\leq p\leq s}\kappa_{s,r}$.

Recall that for $\gamma\in\BN[I]$ we denote by $\Gamma(\gamma)$ the set
of all partitions of $\gamma$, i.e. multisubsets (subsets with multiplicities)
$\Gamma=\lbr \gamma_1,\ldots,\gamma_k\rbr $ of $\BN[I]$ with
$\sum_{r=1}^k\gamma_r=\gamma,\ \gamma_r>0$ (see e.g. ~\cite{k}, 1.3).

The configuration space of colored effective divisors of
multidegree $\gamma$ (the set of colors is $I$) is denoted by $C^\gamma$.
The diagonal stratification $C^\gamma=\sqcup_{\Gamma\in\Gamma(\gamma)}
C^\gamma_\Gamma$ was introduced e.g. in {\em loc. cit.} Recall that for
$\Gamma=\lbr \gamma_1,\ldots,\gamma_k\rbr $ we have $\dim C^\gamma_\Gamma=k$.

\subsubsection{}
For the definition of Laumon's Quasiflags' space $\CQ^L_\alpha$ the reader
may consult ~\cite{la} 4.2, or ~\cite{k} 1.4. It is the space of complete
flags of locally free subsheaves
$$0\subset E_1\subset\dots\subset E_{n-1}\subset V\otimes\CO_C=:\CV$$
such that rank$(E_k)=k$, and $\deg(E_k)=-a_k$.

It is known to be a smooth
projective variety of dimension $2|\alpha|+\dim\CB$.

\subsubsection{}
For the definition of Drinfeld's Quasimaps' space $\CQ^D_\alpha$ the
reader may consult ~\cite{k} 1.2. It is the space of collections of
invertible subsheaves $\CL_\lambda\subset V_\lambda\otimes\CO_C$ for
each dominant weight $\lambda\in X^+$ satisfying Pl\"ucker relations,
and such that $\deg\CL_\lambda=-\langle\lambda,\alpha\rangle$.

It is known to be a (singular, in general) projective variety of
dimension $2|\alpha|+\dim\CB$.

The open subspace $\CQ_\alpha\subset\CQ^D_\alpha$ of genuine maps is formed
by the collections of line subbundles (as opposed to invertible subsheaves)
$\CL_\lambda\subset V_\lambda\otimes\CO_C$. In fact, it is an open stratum
of the stratification by the {\em type of degeneration} of $\CQ^D_\alpha$
introduced in ~\cite{k} 1.3:
$$\CQ^D_\alpha=\bigsqcup_{\beta\leq\alpha}^{\Gamma\in\Gamma(\alpha-\beta)}
\fD^{\beta,\Gamma}_\alpha$$
We have $\fD_{\alpha,\emptyset}=\CQ_\alpha$, and $\fD^{\beta,\Gamma}_\alpha=
\CQ_\beta\times C^{\alpha-\beta}_\Gamma$ (see {\em loc. cit.} 1.3.5).

The space $\CQ^D_\alpha$ is naturally embedded into the product of projective
spaces
$$\BP_\alpha=\prod_{1\leq p\leq n-1}
\BP(\Hom(\CO_C(-\langle\omega_p,\alpha\rangle),
V_{\omega_p}\otimes\CO_C))$$ and is closed in it (see {\em loc. cit.} 1.2.5).
Here $\omega_p$ stands for the fundamental weight dual to the coroot $i_p$.
The fundamental representation $V_{\omega_p}$ equals $\Lambda^pV$.

\subsection{}
\label{main}
We will study the characteristic cycle of the Goresky-MacPherson
perverse sheaf (or the corresponding regular holonomic $D$-module) $IC_\alpha$
on $\CQ^D_\alpha$. As $\CQ^D_\alpha$ is embedded into the smooth space
$\BP_\alpha$, we will view this characteristic cycle $SS(IC_\alpha)$ as a
Lagrangian cycle in the cotangent bundle $T^*\BP_\alpha$. {\em A priori} we
have the following equality:
$$SS(IC_\alpha)=\overline{T^*_{\CQ_\alpha}\BP_\alpha}+
\sum_{\beta<\alpha}^{\Gamma\in\Gamma(\alpha-\beta)}m^{\beta,\Gamma}_\alpha
\overline{T^*_{\fD^{\beta,\Gamma}_\alpha}\BP_\alpha},$$
closures of conormal bundles with multiplicities.

{\bf Theorem.}
$SS(IC_\alpha)=\overline{T^*_{\CQ_\alpha}\BP_\alpha}$ is irreducible.

In the following subsections we will reduce the Theorem to a statement about
geometry of Laumon's resolution.

\subsection{}
We fix a coordinate $z$ on $C$ identifying it with the standard $\BP^1$.
We denote by $\CQ^\infty_\alpha\subset\CQ^D_\alpha$ the open subspace formed by
quasimaps which are genuine maps in a neighbourhood of the point $\infty\in C$.
In other words, $(\CL_\lambda\subset V_\lambda\otimes\CO_C)_{\lambda\in X^+}
\in\CQ^\infty_\alpha$ iff for each $\lambda$ the invertible subsheaf
$\CL_\lambda\subset V_\lambda\otimes\CO_C$ is a line subbundle in some
neighbourhood of $\infty\in C$.

Evidently, $\CQ^\infty_\alpha$ intersects all the strata $\fD_{\beta,\Gamma}$.
Thus it suffices to prove the irreducibility of the singular support of
Goresky-MacPherson sheaf of $\CQ^\infty_\alpha$.

There is a well-defined map of evaluation at $\infty\in C$:
$$\Upsilon_\alpha:\ \CQ^\infty_\alpha\lra\CB$$
It is compatible with the stratification of $\CQ^\infty_\alpha$ and realizes
$\CQ^\infty_\alpha$ as a (stratified) fibre bundle over $\CB$.
In effect, $G$ acts naturally both on $\CQ^\infty_\alpha$
(preserving stratification) and on $\CB$;
the map $\Upsilon_\alpha$ is equivariant, and $\CB$ is homogeneous. We denote
the fiber $\Upsilon_\alpha^{-1}(B)$ over the point $B\in\CB$ by $\CZ_\alpha$.

It inherits the stratification $$\CZ_\alpha=
\bigsqcup_{\beta\leq\alpha}^{\Gamma\in\Gamma(\alpha-\beta)}
\CZ\fD^{\beta,\Gamma}_\alpha$$ from $\CQ^\infty_\alpha$ and $\CQ^D_\alpha$.
It is just the transversal intersection of the fiber $\Upsilon_\alpha^{-1}(B)$
with the stratification of $\CQ^\infty_\alpha$.
As in ~\cite{k} 1.3.5 we have $\CZ\fD^{\beta,\Gamma}_\alpha\iso\CZ_\beta\times
(C-\infty)^{\alpha-\beta}_\Gamma$.

Hence it suffices to prove the irreducibility of the singular support
$SS(IC(\CZ_\alpha))$ of Goresky-MacPherson sheaf $IC(\CZ_\alpha)$ of
$\CZ_\alpha$.

\subsection{Factorization}
The Theorem 6.3 of ~\cite{fm} admits the following immediate Corollary.
Let $(\phi_\beta,\gamma_1x_1,\ldots,\gamma_kx_k)=\phi_\alpha\in
\CZ_\beta\times(C-\infty)^{\alpha-\beta}_\Gamma=\CZ\fD^{\beta,\Gamma}_\alpha
\subset\CZ_\alpha$. Consider also the points $(\phi_r,\gamma_rx_r)=
\phi_{\gamma_r}\in\CZ_0\times(C-\infty)^{\gamma_r}_{\{\{\gamma_r\}\}}=
\CZ\fD^{0,\{\{\gamma_r\}\}}_{\gamma_r}\subset\CZ_{\gamma_r},\ 1\leq r\leq k$.

{\bf Proposition.} There is an analytic open neighbourhood $U_\alpha$
(resp. $U_\beta$, resp. $U_{\gamma_r},\ 1\leq r\leq k$) of $\phi_\alpha$
(resp. $\phi_\beta$, resp. $\phi_{\gamma_r},\ 1\leq r\leq k$) in $\CZ_\alpha$
(resp. $\CZ_\beta$, resp. $\CZ_{\gamma_r},\ 1\leq r\leq k$) such that
$$U_\alpha\iso U_\beta\times\prod_{1\leq r\leq k}U_{\gamma_r}$$
$\Box$

Recall the nonnegative integers $m^{\beta,\Gamma}_\alpha$ introduced in
~\ref{main}. The Proposition implies the following Corollary.

{\bf Corollary.} $m^{\beta,\Gamma}_\alpha=\prod_{1\leq r\leq k}
m^{0,\{\{\gamma_r\}\}}_{\gamma_r}$. $\Box$

Thus to prove that all the multiplicities $m^{\beta,\Gamma}_\alpha$ vanish,
it suffices to check the vanishing of $m^{0,\{\{\gamma\}\}}_\gamma$ for
arbitrary $\gamma>0$.

\subsection{}
It remains to prove that the conormal bundle
$T^*_{\fD^{0,\{\{\gamma\}\}}_\gamma}\BP_\alpha$ to the closed stratum of
$\CQ_\gamma$ enters the singular support $SS(IC_\alpha)$ with multiplicity 0.
To this end we choose a point $(B,\gamma0)=\phi\in\CB\times C=
\CQ_0\times C^\gamma_{\{\{\gamma\}\}}=\fD_\gamma^{0,\{\{\gamma\}\}}
\subset\CQ_\gamma\subset\BP_\gamma$. We also choose a sufficiently generic
meromorphic function $f$ on $\BP_\gamma$ regular around $\phi$ and vanishing
on $\fD_\gamma^{0,\{\{\gamma\}\}}$. According to the Proposition 8.6.4 of
~\cite{ks}, the multiplicity in question is 0 iff $\Phi_f(IC_\gamma)_\phi=0$,
i.e. the stalk of vanishing cycles sheaf at the point $\phi$ vanishes.

To compute the stalk of vanishing cycles sheaf we use the following argument,
borrowed from ~\cite{bfl} ~\S1. As $\pi:\ \CQ^L_\gamma\lra\CQ^D_\gamma$ is a
small resolution of singularities, up to a shift, $IC_\alpha=\pi_*\uQ$.
By the proper base change, $\Phi_f\pi_*\uQ=\pi_*\Phi_{f\circ\pi}\uQ$.
So it suffices to check that $\Phi_{f\circ\pi}\uQ|_{\pi^{-1}(\phi)}=0$.

Let us denote the differential of the function $f$ at the point $\phi$ by
$\xi$ so that $(\phi,\xi)\in T^*_{\fD_\gamma^{0,\{\{\gamma\}\}}}\BP_\gamma$.
Then the support of $\Phi_{f\circ\pi}\uQ|_{\pi^{-1}(\phi)}$ is {\em a priori}
contained in the {\em microlocal fiber} over $(\phi,\xi)$ which we define
presently.

\subsubsection{Definition} Let $\varpi:\ A\to B$ be a map of smooth varieties.
For $a\in A$ let $d_a^*\varpi:\ T^*_{\varpi(a)}B\lra T^*_aA$
denote the codifferential,
and let $(b,\eta)$ be a point in $T^*B$. Then the {\em microlocal fiber}
of $\varpi$ over $(b,\eta)$ is defined to be the set of points
$a\in \varpi^{-1}(b)$
such that $d^*_a\varpi(\eta)=0$.

\subsubsection{}
\label{prop}
Thus we have reduced the Theorem ~\ref{main} to the following Proposition.

{\bf Proposition.} For a sufficiently generic $\xi$ such that $(\phi,\xi)
\in T^*_{\fD_\gamma^{0,\{\{\gamma\}\}}}\BP_\gamma$, the microlocal fiber
of Laumon's resolution $\pi$ over $(\phi,\xi)$ is empty. Equivalently,
the cone $\cup_{E_\bullet\in\pi^{-1}(\phi)}\Ker (d^*_{E_\bullet}\pi)$ is a
proper subvariety of the fiber of
$T^*_{\fD_\gamma^{0,\{\{\gamma\}\}}}\BP_\gamma$ at $\phi$.

\subsection{Piecification of a simple fiber}
The fiber $\pi^{-1}(\phi)$ was called the {\em simple fiber} in ~\cite{k} ~\S2.
It was proved in {\em loc. cit.} ~2.3.3 that $\pi^{-1}(\phi)$ is a disjoint
union of (pseudo)affine spaces $\fS(\mu(\kappa))$ where $\kappa$ runs through
the set $\fK(\gamma)$ of Kostant partitions of $\gamma$ (for the notation
$\mu(\kappa)$ see ~\ref{not} or ~\cite{k} ~(9)). Another way to parametrize
these pseudoaffine pieces was introduced in ~\cite{fk} ~2.11. Let us recall
it here.

We define nonnegative integers $c_p, 1\leq p\leq n-1$, so that
$\gamma=\sum_{p=1}^{n-1}c_pi_p$.

\subsubsection{Definition} $\CalD(\gamma)$ is the set of collections of
nonnegative integers $(d_{p,q})_{1\leq q\leq p\leq n-1}$ such that

a) For any $1\leq q\leq p\leq r\leq n-1$ we have $d_{r,q}\leq d_{p,q}$;

b) For any $1\leq p\leq n-1$ we have $\sum_{q=1}^pd_{p,q}=c_p$.

\subsubsection{Lemma} The correspondence $\kappa=(\kappa_{p,q})_{1\leq q\leq p
\leq n-1}\mapsto(d_{p,q}:=\sum_{r=p}^{n-1}\kappa_{r,q})_{1\leq q\leq p
\leq n-1}$ defines a bijection between $\fK(\gamma)$ and $\CalD(\gamma)$. $\Box$

\subsubsection{} Using the above Lemma we can rewrite the parametrization
of the pseudoaffine pieces of the simple fiber as follows:
$$\pi^{-1}(\phi)=\bigsqcup_{\fd\in\CalD(\gamma)}\fS(\fd)$$
In these terms the dimension formula of ~\cite{k} ~2.3.3 reads as follows:
for $\fd=(d_{p,q})_{1\leq q\leq p\leq n-1}$ we have $\dim\fS(\fd)=
\sum_{1\leq q<p\leq n-1}d_{p,q}$.

Note also that $\sum_{1\leq q\leq p\leq n-1}d_{p,q}=
\sum_{1\leq p\leq n-1}c_p=|\gamma|$.

\subsection{Proposition}
\label{red}
For arbitrary $\fd=(d_{p,q})_{1\leq q\leq p\leq n-1}\in\CalD(\gamma)$ and
arbitrary quasiflag $E_\bullet\in\fS(\fd)\subset\pi^{-1}(\phi)$ we have
$\dim\Ker (d_{E_\bullet}\pi)<\sum_{1\leq p\leq n-1}d_{p,q}
+\sum_{1\leq q\leq p\leq n-1}d_{p,q}-1$.

This Proposition implies the Proposition ~\ref{prop} straightforwardly.
In effect, $\codim\Ker (d_{E_\bullet}\pi)=\dim\CQ^L_\gamma-\dim\Ker
(d_{E_\bullet}\pi)>2|\gamma|+\dim\CB-\sum_{1\leq p\leq n-1}d_{p,p}
-\sum_{1\leq q\leq p\leq n-1}d_{p,q}+1=
\dim\CB+1+\sum_{1\leq q<p\leq n-1}d_{p,q}$.
Hence the subspace $\Ker (d^*_{E_\bullet}\pi)\subset T^*_\phi\BP_\gamma$ has
codimension greater than $\dim\CB+1+\sum_{1\leq q<p\leq n-1}d_{p,q}$.
Recall that $\dim\fD_\gamma^{0,\{\{\gamma\}\}}=\dim\CB+1$.
Hence the codimension of $\Ker (d^*_{E_\bullet}\pi)\cap
T^*_{\fD_\gamma^{0,\{\{\gamma\}\}}}\BP_\gamma$ in the fiber of
$T^*_{\fD_\gamma^{0,\{\{\gamma\}\}}}\BP_\gamma$ at $\phi$ is greater than
$\sum_{1\leq q<p\leq n-1}d_{p,q}=\dim\fS(\fd)$. Hence
the cone $\cup_{E_\bullet\in\fS(\fd)}\Ker (d^*_{E_\bullet}\pi)$ is a
proper subvariety of the fiber of
$T^*_{\fD_\gamma^{0,\{\{\gamma\}\}}}\BP_\gamma$ at $\phi$.

The union of these proper subvarieties over $\fd\in\CalD(\gamma)$ is again a
proper subvariety of the fiber of
$T^*_{\fD_\gamma^{0,\{\{\gamma\}\}}}\BP_\gamma$ at $\phi$ which concludes the
proof of the Proposition ~\ref{prop}.

\subsection{Fixed points}
It remains to prove the Proposition ~\ref{red}. To this end recall that
the Cartan group $H$ acts on $V$ and hence on $\CQ^L_\alpha$.
The group $\BC^*$ of dilations of $C={\Bbb P}^1$ preserving $0$ and $\infty$
also acts on $\CQ^L_\alpha$ commuting with the action of $H$.
Hence we obtain the action of a torus $\BT:=H\times\BC^*$ on $\CQ^L_\alpha$.

It preserves the simple fiber $\pi^{-1}(\phi)$ and its pseudoaffine pieces
$\fS(\fd),\ \fd\in\CalD(\gamma)$, for evident reasons. It was proved in
~\cite{fk} ~2.12 that each piece
$\fS(\fd),\ \fd=(d_{p,q})_{1\leq q\leq p\leq n-1}$ contains exactly
one $\BT$-fixed point $\delta(\fd)=(E_1,\ldots,E_{n-1})$. Here
$$
{\arraycolsep=1pt
\begin{array}{llrlcrlcccrlc}
E_1 & = E_{1,1} \\
E_2 & = E_{2,1} &\opl& E_{2,2} \\
\ \vdots && \vdots &&& \vdots \\
E_{n-1} & = E_{n-1,1} &\opl& E_{n-1,2} &\opl& \dots &\opl& E_{n-1,n-1} \\
\end{array}
}
$$
and $E_{p,q}=\CO(-d_{p,q})\subset\CO v_q\subset\CV=V\otimes\CO_C$
with quotient sheaf $\dfrac\CO{\CO(-d_{p,q})}$ concentrated at $0\in C$.

%
%

\subsubsection{}
\label{key}
Now the $\BT$-action contracts $\fS(\fd)$ to $\delta(\fd)$. Since the map
$\pi$ is $\BT$-equivariant, and the dimension of $\Ker (d_{E_\bullet}\pi)$
is lower semicontinuous, the Proposition ~\ref{red} follows from the next one.

{\bf Key Proposition.}
For arbitrary $\fd=(d_{p,q})_{1\leq q\leq p\leq n-1}\in\CalD(\gamma)$
($\gamma\ne0$) we have
$\dim\Ker (d_{\delta(\fd)}\pi)<\sum_{1\leq p\leq n-1}d_{p,p}
+\sum_{1\leq q\leq p\leq n-1}d_{p,q}-1$.

The proof will be given in the next section.

\subsubsection{Remark} In general, the pieces $\fS(\fd)$ of the simple fiber
are not equisingular, i.e. $\dim\Ker (d_{E_\bullet}\pi)$ is not constant along
a piece. The simplest example occurs for $G=SL_3,\ \gamma=2i_1+2i_2$.
Then the simple fiber is a singular 2-dimensional quadric. Its singular point
is the fixed point of the 1-dimensional piece $\fS(\fd)$ where
$d_{1,1}=2,d_{2,1}=d_{2,2}=1$.
At this point we have $\dim\Ker (d_{\delta(\fd)}\pi)=3$ while at the other
points in this piece we have $\dim\Ker (d_{E_\bullet}\pi)=2$.

\section{The proof of the Key Proposition}

\subsection{Tangent spaces}
Let $\Omega$ be the following quiver: $\Omega=1\lra 2\lra\ldots\lra n-1$.
Thus the set of vertices coincides with $I$. A quasiflag $(E_1\hra E_2\hra
\ldots\hra E_{n-1}\subset\CV)\in\CQ^L_\gamma$ may be viewed as a representation
of $\Omega$ in the category of coherent sheaves on $C$. If we denote the
quotient sheaf $\CV/E_p$ by $Q_p,\ 1\leq p\leq n-1$, we have another
representation of $\Omega$ in coherent sheaves on $C$, namely,
$$Q_\bullet:=(Q_1\twoheadrightarrow Q_2\twoheadrightarrow
\ldots\twoheadrightarrow Q_{n-1})$$

\subsubsection{Exercise} $T_{E_\bullet}\CQ^L_\gamma=\Hom_\Omega(E_\bullet,
Q_\bullet)$ where $\Hom_\Omega(?,?)$ stands for the morphisms in the category
of representations of $\Omega$ in coherent sheaves on $C$.

\subsubsection{} Consider a point $\CL_\bullet=
(\CL_1,\ldots,\CL_{n-1})\in\BP_\gamma$.
Here $\CL_p\subset V_{\omega_p}\otimes\CO_C$ is an invertible subsheaf,
the image of morphism $\CO_C(-\langle\omega_p,\gamma\rangle)\hra V_{\omega_p}
\otimes\CO_C$.

{\em Exercise.} $T_{\CL_\bullet}\BP_\gamma=\prod_{p=1}^{n-1}\Hom(\CL_p,
V_{\omega_p}\otimes\CO_C/\CL_p)$.

\subsubsection{} Recall that for $E_\bullet\in\CQ^L_\gamma$ we have
$\pi(E_\bullet)=\CL_\bullet\in\BP_\gamma$ where
$\CL_p=\Lambda^pE_p$ for $1\leq p\leq n-1$.

{\em Exercise.} For $h_\bullet=(h_1,\ldots,h_{n-1})\in
T_{E_\bullet}\CQ^L_\gamma$ we have $d_{E_\bullet}\pi(h_\bullet)=
(\Lambda^1h_1,\Lambda^2h_2,\ldots,\Lambda^{n-1}h_{n-1})\in
T_{\CL_\bullet}\BP_\gamma$.

\subsection{}
From now on we fix $\gamma>0,\ \fd\in\CalD(\gamma),\ \delta(\fd)=:E_\bullet$.
To unburden the notations we will denote the tangent space
$T_{E_\bullet}\CQ^L_\gamma$ by $T$.
Since $\QL\gamma$ is a smooth $(2|\gamma|+\dim\CB)$-dimensional variety
it suffices to find a subspace $N\subset T$ of dimension
$$
2|\gamma|+\dim\CB-\sum_{1\le p\le n-1}d_{p,p}-\sum_{1\le q\le p\le n-1}d_{p,q}+1=
\sum_{1\le q<p\le n-1}(d_{p,q}+1)+1
$$
such that $d_{E_\bullet}\pi|_{N}$ is injective.

\subsection{}
Let $N_0=\opll_{n-1\ge p>q\ge 1}\Hom(\CO(-d_{p,q}),\CO)$.
We have $\dim N_0=\sum_{n-1\ge p>q\ge 1}(d_{p,q}+1)$.

Recall that we have canonically
$T=\Hom_\Om(E_\bullet,Q_\bullet)$, where
$$
Q_p=\CV/E_p=\left(\opll_{q=1}^p\left(\dfrac\CO{\CO(-d_{p,q})}\right)v_q\right)
\opl\left(\opll_{q=p+1}^n\CO v_q\right).
$$

\subsection{}
Let us define a map $\nu_0:N_0\to T$ assigning to an element
$(f_{p,q})\in N_0$ a morphism $\nu_0(f_{p,q}):=
F\in\Hom_\bullet(E_\bullet,Q_\bullet)$ of graded coherent sheaves, where
$F|_{E_{p,q}}=\opll_{r=p+1}^n F_{p,q}^r$, and
$$
F_{p,q}^r:E_{p,q}\to\CO v_r\subset Q_p\quad
\text{is defined as the composition}\quad
E_{p,q} \subset E_{r,q}=\CO(-d_{r,q}) @>f_{r,q}>> \CO v_r
$$

\subsection{}{\bf Lemma.}
The map $F:E_\bullet\to Q_\bullet$ is a morphism
of representations of the quiver $\Om$.
\begin{pf}
We need to check the commutativity of the following diagram
$$
\begin{CD}
E_p @>>>    E_{p'}  \\
@VFVV       @VFVV   \\
Q_p @>>>    Q_{p'}
\end{CD}
$$
Since $E_p$ and $Q_{p'}$ are canonically decomposed into the direct
sum it suffices to note that for any $q\le p\le p'<r$ the following diagram
$$
\begin{CD}
E_{p,q}     @>>>    E_{p',q}    @>>>    E_{r,q}     \\
@VF_{p,q}^rVV       @VF_{p',q}^rVV      @Vf_{r,q}VV \\
\CO v_r     @=  \CO v_r     @=  \CO v_r
\end{CD}
$$
commutes and for any $q\le p<r\le p'$ the following diagram
$$
\begin{CD}
E_{p,q}     @>>>    E_{r,q}     @>>>    E_{p',q}    \\
@VF_{p,q}^rVV       @Vf_{r,q}VV     @V0VV       \\
\CO v_r     @=  \CO v_r     @>>>
\left(\dfrac\CO{\CO(-d_{r,q})}\right)v_r
\end{CD}
$$
commutes as well.
\end{pf}

\subsection{}

Let $N_1={\Bbb C}$. Let $p_0=\min\{1\le p\le n-1\ |\ d_{p,p}>0\}$
and pick a non-zero element $\ff\in\Hom(\CO(-d_{p,p}),\frac\CO{\CO(-d_{p,p})})$.
Define the map $\nu_1:N_1\to T$ by assigning to $1\in N_1$
the element $\fF\in\Hom_\Om(E_\bullet,Q_\bullet)$ defined on
$E_{p,p}$ as the composition
$$
E_{p,p}=\CO(-d_{p,p}) @>\ff>> \frac\CO{\CO(-d_{p,p})}v_p\subset Q_p
$$
and with all other components equal to zero.

\subsection{}
Let $\CM(r,d;\CV)$ denote the space of
rank $r$ and degree $d$ subsheaves in $\CV$.

\subsubsection{}
\label{lem}
Let $E\subset\CV$ be a rank $k$ and degree $d$ subsheaf in the
vector bundle $\CV$. Let $\CV/E=\CT\opl \CF$ be a decomposition of
the quotient sheaf into the sum of the torsion $\CT$ and a locally
free sheaf $\CF$.
Consider the map $\det:\CM(r,d;\CV)\to\CM(1,d;\Lambda^k\CV)$
sending $E$ to $\Lambda^kE$. Then the restriction of its differential
$d_E\det:T_E\CM(r,d;\CV)=\Hom(E,\CV/E)\to
\Hom(\Lambda^kE,\Lambda^kV/\Lambda^kE)=T_{\Lambda^kE}\CM(1,d;\Lambda^k\CV)$
to the subspace $\Hom(E,\CF)\subset\Hom(E,\CV/E)$ factors as
$\Hom(E,\CF)\cong\Hom(\Lambda^kE,\Lambda^{k-1}E\ot \CF)\subset
\Hom(\Lambda^kE,\Lambda^kV/\Lambda^kE)$. Therefore it is injective.

\subsubsection{}\label{lem1}
Let $E=\CO^{\opl(r-1)}\opl\CO(-d)$ be a subsheaf in $\CV=\CO^{\opl n}$.
Then the restriction of differential $d_E\det$ to the subspace
$\Hom(E,\CT)\subset\Hom(E,\CV/E)$ is injective.

This immediately follows from the following fact.
Let $\TE=\CO^{\opl r}\subset\CV$ be the normalization of $E$ in $\CV$,
that is, the maximal vector subbundle $\TE\subset\CV$ such that $\TE/E$ is
torsion.
Then $\CT=\TE/E\cong\Lambda^k\TE/\Lambda^kE\subset\Lambda^k\CV/\Lambda^kE$.

\subsubsection{}\label{rem}
Clearly, the subsheaves
$\Lambda^{k-1}E\ot \CF\subset\Lambda^k\CV/\Lambda^kE$
and $\CT\subset\Lambda^k\CV/\Lambda^kE$ do not intersect.

\subsubsection{}
It follows from \ref{lem}, \ref{lem1} and \ref{rem} that the composition
$d_{E_\bullet}\pi\circ(\nu_0\opl\nu_1):N_0\opl N_1\to T_{\pi(E_\bullet)}\CP$
is injective, hence
$N:=(\nu_0\opl\nu_1)(N_0\opl N_1)\subset T_{E_\bullet}\QD\gamma$ enjoys
the desired property. Namely, $d_{E_\bullet}\pi|_N$ is injective, and
$\dim N=\sum_{1\leq q<p\leq n-1}(d_{p,q}+1)+1$.

This completes the proof of the Key Proposition ~\ref{key} along with the
Main Theorem ~\ref{main}. $\Box$

\end{document}